\documentclass[11pt]{amsart}
\usepackage{graphicx,calc,amscd}
\usepackage{epsfig,psfig}
\usepackage{float}
\usepackage{labelfig}
\usepackage{amsmath}
\usepackage{amsfonts}
\usepackage{amssymb}

\newcommand{\R}{{\mathbb R}}

\newcommand{\Hyp}{{\mathbb H}}

\newcommand{\CC}{{\mathcal C}}

\newcommand{\DD}{{\mathbb D}}
\newcommand{\Sp}{{\mathbb S}}

\newcommand{\mm}{\mathrm{min}}

\renewcommand{\tt}{\mbox{\rm Tr}_{2^d|2}\,}

\newcommand{\arccosh}{\mathrm{arccosh}}

\begingroup\makeatletter\ifx\SetFigFont\undefined%
\gdef\SetFigFont#1#2#3#4#5{%
 \reset@font\fontsize{#1}{#2pt}%
 \fontfamily{#3}\fontseries{#4}\fontshape{#5}%
 \selectfont}%
\fi\makeatother\endgroup%
\begin{document}
\title{On angles formed by $N$ points of the Euclidean and Hyperbolic planes}
\author{Ghislain Jaudon and Hugo Parlier}

\address{Section de Math\'ematiques\\
Universit\'e de Gen\`eve\\
1211 Gen\`eve\\
SWITZERLAND} \email{ghislain.jaudon@math.unige.ch}
\email{hugo.parlier@math.unige.ch}

\thanks{The first author was supported in part by the Swiss National Science
Foundation, \\
\indent No.~PP002-68627.}
\date{\today}

\maketitle

\section{Introduction}\label{Sect:S1}

Consider $n$ points in the Euclidean plane, and consider the
angles that are formed by all possible triples of points. In 1939,
L. Blumenthal \cite{bl39} proposed the problem of finding a sharp
lower bound on the greatest of all these angles (which we shall
refer to as the {\it mini-max problem}). The full solution was
given by Bl. Sendov \cite{sen95} more than a half century later,
after several partial resolutions by a variety of authors.
Recently, a very similar problem was proposed in \cite{brmopabook}
(Conjecture 6, p. 446)
where the following is conjectured:\\

{\it Among all angles formed by triples of $n$ points in the
plane there is an angle of at most $\frac{\pi}{n}$.}\\

Furthermore, it is conjectured that the value $\frac{\pi}{n}$ is
only obtained when the $n$ points form a regular $n$-gon. At first
glance, there is no apparent reason why the solution to finding
the maximum least angle formed by $n$ points (which we shall refer
to as the {\it maxi-min problem}) would be any less difficult than
the mini-max problem. The first goal of this note is to show that
the conjecture is correct, and surprisingly the proof turns out to
be elementary.\\

Subsequently, we treat the same question for $n$ points in the
hyperbolic plane. A similar proof shows that any configuration of
$n$ points of the hyperbolic plane form at least one angle
strictly less than $\frac{\pi}{n}$. This value is also shown to be
sharp, but this is slightly more complicated, because there is no
maximal configuration as in the Euclidean case (and thus the
``$\max$"
becomes a ``$\sup$").\\

Before beginning with a proof of the conjecture, let us fix a few
notations which we will carry with us throughout this note. For
example, what we mean by the angle formed by a triple of points is
the following, for both the Euclidean and hyperbolic plane. The
angle associated to the triple $(p,q,r)$, denoted $\angle
(p,q,r)$, is the positive interior angle of the triangle
$\triangle (p,q,r)$ at point $q$. By our definition $\angle
(p,q,r)=\angle (r,q,p)\in [0,\pi]$ and $\angle (p,q,r)\in ]0,\pi[$
if $\triangle (p,q,r)$ is not degenerated. Let $M$ be some metric
space where the notions of triple of points and angle formed by a
triple make some kind of sense. We shall denote a configuration of
$n$ points in $M$ by $\CC$ and by $\alpha_\mm(\CC)$ the smallest
angle formed by all triples of $\CC$. We further denote by
$\CC_n(M)$ the set of all possible configurations of $n$ points.
Finally, $\alpha_n(M):=\sup_{\CC\in\CC_n(M)}\alpha_\mm(\CC)$. For
instance, in the following section we shall show that
$\alpha_n(\R^2)=\frac{\pi}{n}$ and that the ``sup" is in fact a
``max" attained by a unique type of configuration.

\section{Minimum angles in the plane}\label{Sect:S2}
Our first goal is to show that $n$ distinct points in the plane
always form at least one angle which is at most $\frac{\pi}{n}$.
Let us suppose that this is not the case, i.e., there exists a
configuration $\mathcal{C}$ of $n$ points such that all the
$3\binom{n}{3}$ angles formed by the $\binom{n}{3}$ triples of
points are all strictly superior to $\frac{\pi}{n}$. In
particular, every angle formed by three points in $\mathcal{C}$ is
in $]0,\pi[$.
\begin{figure}[h]
\leavevmode \SetLabels
\L(.5*.06) $\gamma_1$\\
\L(.44*.11) $\gamma_2$\\
\L(.517*.45) $\gamma_3$\\
\L(.43*.38) $\gamma_4$\\
\L(.37*.4) $\gamma_5$\\
\L(.325*.51) $\gamma_6$\\
\L(.6*.03) $x_1$\\
\L(.73*.26) $x_2$\\
\L(.62*.614) $x_3$\\
\L(.66*.965) $x_4$\\
\L(.43*.79) $x_5$\\
\L(.32*1.02) $x_6$\\
\L(.26*.8) $x_7$\\
\L(.37*-.04) $p$\\
\endSetLabels
\begin{center}
\AffixLabels{\centerline{\epsfig{file =
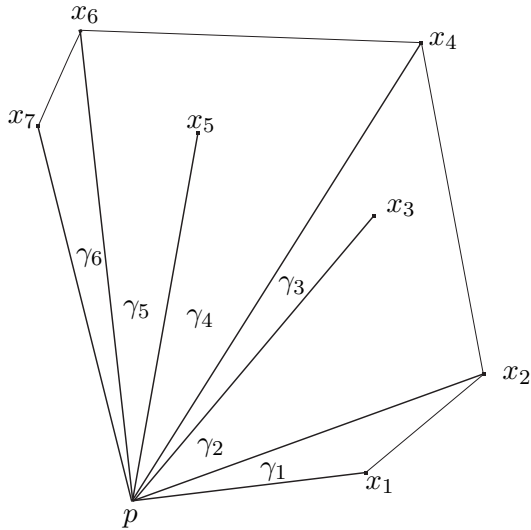,width=6cm,angle=0}}}
\end{center}
\caption{The convex hull of $8$ points and the angles $\gamma_i$}
\label{fig:convexhull}
\end{figure}

Take an extremal point $p$ in the convex hull $X$ of
$\mathcal{C}$. Let us denote by $x_1$ and $x_{n-1}$ the two
extremal points adjacent to $p$ on the boundary of $X$ (by
hypothesis, $X$ has at least three extremal points). Consider an
enumeration $x_2,\ldots, x_{n-2}$ of the remaining points of
$\mathcal{C}$ such that $\angle (x_1,p,x_2)< \angle (x_1,p,x_3)<
\ldots < \angle (x_1,p,x_{n-1})$. We denote $\gamma_i:=\angle
(x_i,p,x_{i+1})$, for $i=1,\hdots, n-2$. Now, according to our
assumption on $\mathcal{C}$, each angle $\gamma_i$  is strictly
greater that $\frac{\pi}{n}$. Therefore
\begin{equation}\label{eqn:totalangle}
\angle (x_1,p,x_{n-1}) =\sum_{i=1}^{n-2}\gamma_i>
\frac{(n-2)\pi}{n}.
\end{equation}

However we have assumed that $\angle
(x_1,x_{n-1},p)>\frac{\pi}{n}$ and $\angle
(p,x_1,x_{n-1})>\frac{\pi}{n}$, and thus we get the following
contradiction:
$$
\pi=\angle (x_1,x_{n-1},p) + \angle (x_1,p,x_{n-1}) + \angle
(p,x_1,x_{n-1}) >
\frac{\pi}{n}+\frac{(n-2)\pi}{n}+\frac{\pi}{n}=\pi.
$$
\\
This shows that any configuration of $n$ points in the plane
necessarily contains an angle smaller than or equal to
$\frac{\pi}{n}$. Moreover, if $\mathcal(C)$ is a regular $n$-gon,
all points on the boundary of $X$, and for any choice of point
$p$, the angles $\gamma_i$ are all equal to $\frac{\pi}{n}$. The
smallest angle is exactly equal to $\frac{\pi}{n}$ which shows
that the bound is sharp. The next step is to show that if all of
the angles are greater than or equal to $\frac{\pi}{n}$, then this
is the only possible configuration, i.e., the $n$ points are the
vertices of a regular $n$-gon.\\

To show this, let us backtrack a little supposing this time that
$\alpha_\mm(\CC)=\frac{\pi}{n}$. We can consider the convex hull
$X$ of $\CC$, an extremal point $p$ of $X$, and the same
enumeration of the $n-1$ remaining points. In this case, equation
(\ref{eqn:totalangle}) becomes $\angle (x_1,p,x_{n-1})
=\sum_{i=1}^{n-2}\gamma_i\geq \frac{(n-2)\pi}{n}$. By the same
reasoning, if $\angle (x_1,p,x_{n-1})>\frac{(n-2)\pi}{n}$, then
either angle $\angle (x_1,x_{n-1},p)$ or angle $\angle
(p,x_1,x_{n-1})$ is strictly less than $\frac{\pi}{n}$, a
contradiction. Thus $\angle (x_1,p,x_{n-1})=\frac{(n-2)\pi}{n}$,
$\angle (x_1,x_{n-1},p)=\angle (p,x_1,x_{n-1})=\frac{\pi}{n}$ and
the triangle $\triangle(p,x_1,x_{n-1})$ is isosceles. This
reasoning can be applied to every extremal point on the convex
hull, i.e., every extremal point $q$ of $X$ forms an isosceles
triangle with its neighbors of the convex hull with interior angle
in $q$ equal to $\frac{(n-2)\pi}{n}$. It follows that the extremal
points of the convex hull form a regular $N$-gon, with interior
angles equal to $\frac{(n-2)\pi}{n}$. As the interior angle of a
regular $N$-gon is exactly equal to $\frac{(N-2)\pi}{N}$, it
follows that $N=n$, and thus the points of $\CC$ form a regular
$n$-gon. This proves the conjecture.

\section{Angles of points in the hyperbolic plane}

Interestingly, the solution to the problem in the hyperbolic plane
$\Hyp$ is very similar to the problem in the Euclidean plane. Our
goal is not to fully present this complicated and beautiful metric
space, but for clarity, let us recall a few facts concerning the
hyperbolic plane. One way of viewing the hyperbolic plane is as
the interior of the unit disk $\DD=\{x\in \R^2 \mid \|x\|<1\}$
endowed with a metric of constant curvature $-1$. The traces of
its geodesics are the intersection with $\DD$ of the circles and
lines of $\R^2$ which are perpendicular to the boundary of $\DD$
(i.e., $\partial\DD=S^1$) each time they intersect it. In
particular, these circles and lines cross $\partial\DD$ exactly
twice, and the lines pass through the origin. As in $\R^2$, there
is a unique geodesic between two distinct points, and thus a
unique way of defining the triangle associated to a triple of
points $(a,b,c)$. Hyperbolic triangles differ from Euclidean ones
in that the sum of their interior angles, say $\angle_a$,
$\angle_b$ and $\angle_c$, verifies

\begin{equation}\label{eqn:hyptri}
\angle_a+\angle_b+\angle_c<\pi.
\end{equation}

Conversely, for any triple of real positive numbers
$(\theta_a,\theta_b,\theta_c)$ whose sum is strictly less than
$\pi$, there exists a hyperbolic triangle whose interior angles
are exactly equal to $(\theta_a,\theta_b,\theta_c)$ and this
triangle is defined uniquely up to isometry. Furthermore, the area
of a hyperbolic triangle $\triangle$ of interior angles
$(\theta_a,\theta_b,\theta_c)$ is given by

\begin{equation}\label{eqn:area}
{\rm area}(\triangle)=\pi-(\theta_a+\theta_b+\theta_c).
\end{equation}

Essentially, the smaller the interior of hyperbolic triangle gets,
the more Euclidean it looks. For simplicity, our proof will only
make use of these facts.\\

Our goal is to prove that $\alpha_n(\Hyp)=\frac{\pi}{n}$. Let us
begin by showing that for any given configuration $\CC$,
$\alpha_\mm(\CC)<\frac{\pi}{n}$. To show this, suppose that there
is a configuration $\CC$ with $\alpha_\mm(\CC)\geq \frac{\pi}{n}$,
and let us adopt the same notations as in the Euclidean case: the
convex hull of our points will be denoted $X$, an extremal point
of the convex hull $p$ and the rest of the points of $\CC$ denoted
$x_i$, $i=1,\hdots,n-1$. (In fact, figure \ref{fig:convexhull}
remains valid up, albeit that the geodesics shouldn't be so
straight.) We can now mimick the proof in the Euclidean case, and
we obtain that $\angle (x_1,p,x_{n-1})\geq \frac{(n-2)\pi}{n}$.
Now considering the triangle $\triangle(p,x_1,x_{n-1})$ and by
equation (\ref{eqn:hyptri}), $\min\{\angle (x_1,x_{n-1},p),\angle
(p,x_1,x_{n-1}) \}<\frac{\pi}{n}$.\\

\begin{figure}[h]
\leavevmode \SetLabels
\L(.695*.353) $\frac{\pi}{n}$\\
\L(.662*.087) $\frac{\theta_n}{2}$\\
\L(.72*.087) $\frac{\theta_n}{2}$\\
\L(.22*.28) $\theta_a$\\
\L(.11*.49) $\theta_c$\\
\L(.22*.665) $\theta_b$\\
\endSetLabels
\begin{center}
\AffixLabels{\centerline{\epsfig{file =
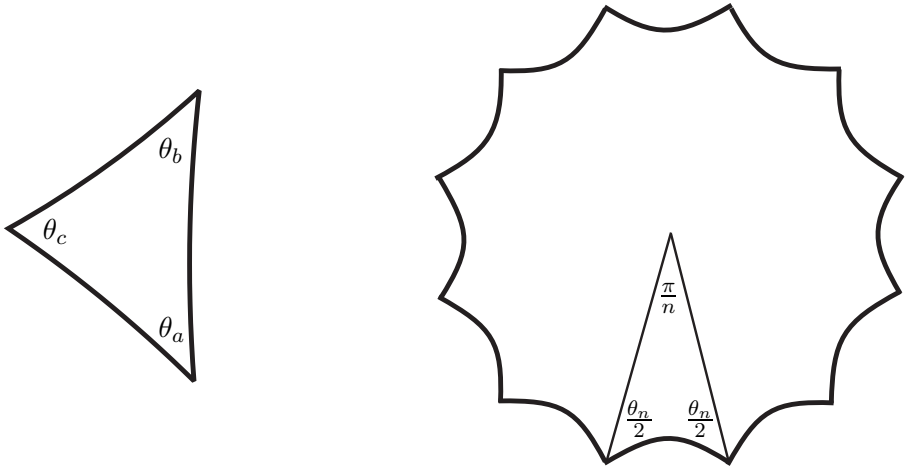,width=12cm,angle=0}}}
\end{center}
\caption{A hyperbolic triangle and a regular $n$-gon with $n=12$}
\label{fig:trianglengon}
\end{figure}

In order to prove $\alpha_n(\Hyp)=\frac{\pi}{n}$, we need to show
that for any $\varepsilon>0$, there exists a configuration
$\CC_\varepsilon$ of $n$ points with
$\angle_\mm(\CC_\varepsilon)\geq \frac{\pi}{n}-\varepsilon$. The
configuration that will prove this is a regular $n$-gon inscribed
in a small disk $D_\varepsilon$ of area $\varepsilon$. (For those
who are interested, the radius $r_\varepsilon$ of such a disk is
equal to $\arccosh\frac{\varepsilon + 2\pi}{2\pi}$.) Any triangle
$\triangle$ of points found inside or on the boundary of the disk
will be entirely contained in the closed disk, and thus
$$
{\rm area}(\triangle)<{\rm area}(D_\varepsilon)=\varepsilon.
$$
By equation (\ref{eqn:area}) it follows that the interior angles
$(\theta_a,\theta_b,\theta_c)$ of $\triangle$ verify

\begin{equation}\label{eqn:littletriangle}
\pi-\varepsilon< \theta_a+\theta_b+\theta_c<\pi.
\end{equation}

Now for $n\geq 3$, consider a regular $n$-gon inscribed in
$D_\varepsilon$, which can be constructed as follows. First chose
a point on the boundary of $D_\varepsilon$, take the unique
geodesic segment between this point and the center of the disk,
and take the image of this segment by rotations of angle
$\frac{2\pi}{n}$ around the center of the disk. The intersection
point of these segments and the boundary of the disks give you the
$n$ vertices of the regular $n$-gon. The polygon, denoted $P_n$,
is then obtained
by taking the $n$ geodesic segments between successive vertices.\\

Consider the triangle formed by two successive vertices and the
center of the disk (see figure \ref{fig:trianglengon}). By what
precedes, the triangle is formed by an edge of the polygon and two
radii of the disk which meet at an angle of $\frac{2\pi}{n}$, and
is thus isosceles. As we know the the length of the two radii,
using hyperbolic trigonometry we could calculate the remaining
lengths and angles of this triangle, but in fact we already know
enough to prove our point. If we denote by $\theta_n$ the interior
angle of our $n$-gon, the two equal interior angles of our
triangle have value $\frac{\theta_n}{2}$. By equation
(\ref{eqn:littletriangle}) applied to this triangle, we have

$$
\frac{n-2}{n}\pi-\varepsilon<\theta_n<\frac{n-2}{n}\pi.
$$

Now we've constructed our $n$-gon, and we have a lower bound on
the interior angles, but what remains to be proved is that any
angle $\gamma$ formed by a triple of vertices of the $n$-gon
verifies $\gamma>\frac{\pi}{n}-\varepsilon$. Let us place
ourselves in a vertex $p$ of $P_n$, and consider our usual
enumeration of our remaining vertices, i.e., $x_1$ and $x_{n-1}$
the two vertices closest to $p$, $x_2,\ldots, x_{n-2}$ of the
remaining vertices of $P_n$ such that $\angle (x_1,p,x_2)< \angle
(x_1,p,x_3)< \ldots < \angle (x_1,p,x_{n-1})=\theta_n$, and denote
by $\gamma_i$ the angle $\angle (x_i,p,x_{i+1})$. First of all, in
the isosceles triangle $\triangle(x_1,p,x_2)$, we have two angles
of value $\gamma_1$, and by equation (\ref{eqn:littletriangle}),
we obtain

$$
2\gamma_1+\theta_n>\pi-\varepsilon,
$$

and thus using what we know about $\theta_n$

$$
\gamma_1>\frac{\pi}{n}-\varepsilon.
$$

If, for a regular $n$-gon in the hyperbolic plane, we had
$\gamma_i=\gamma_1$ for $i=2,\hdots,n-1$ as we do in the Euclidean
plane, the proof would be finished, but in the hyperbolic plane
this is not the case. We could prove this, but we definitely
won't be using it, so we won't.\\

\begin{figure}[h]
\leavevmode \SetLabels
\L(.56*.05) $\theta_n$\\
\L(.49*.074) $\gamma_1$\\
\L(.62*.145) $\gamma_1$\\
\L(.41*-0.03) $p$\\
\L(.58*-0.03) $x_1$\\
\L(.71*.1) $x_2$\\
\L(.73*0.85) $x_k$\\
\L(.578*1.01) $x_{k+1}$\\
\L(.493*.415) $\gamma_k$\\
\L(.655*.816) $\alpha_k$\\
\L(.573*.891) $\beta_k$\\
\L(.28*0.1) $x_{n-1}$\\
\endSetLabels
\begin{center}
\AffixLabels{\centerline{\epsfig{file =
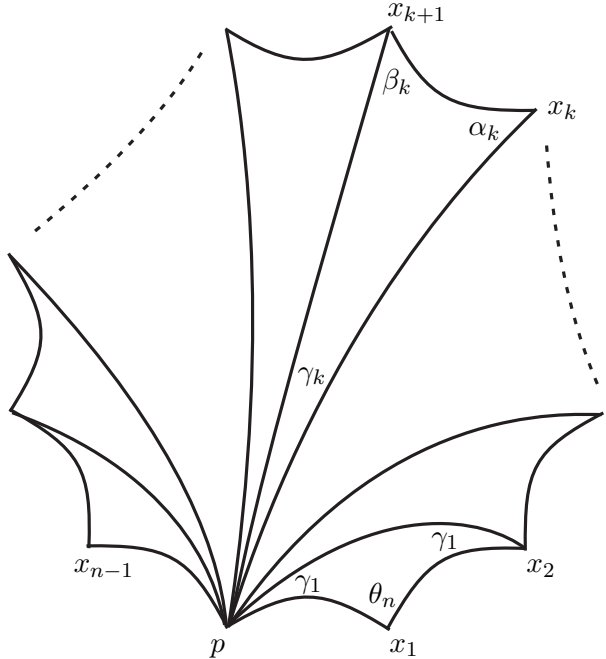,width=8cm,angle=0}}}
\end{center}
\caption{A regular $n$-gon} \label{fig:ngon}
\end{figure}

For $k>1$, consider the triangle $\triangle (x_k,p,x_{k+1})$
which, in virtue of the regularity of $P_n$, has interior angles
of values $\gamma_k$,
$\alpha_k:=\theta_n-\sum_{l=1}^{k-1}\gamma_l$ and
$\beta_k:=\sum_{l=1}^{k}\gamma_l$ (see figure \ref{fig:ngon}). We
can use equation (\ref{eqn:littletriangle}) to obtain

$$
\gamma_k+\sum_{l=1}^{k}\gamma_l +
\theta_n-\sum_{l=1}^{k-1}\gamma_l=2\gamma_k+\theta_n>\pi-\varepsilon
$$

and as before

$$
\gamma_k>\frac{\pi}{n}-\varepsilon.
$$

All other angles formed by triples of the vertices are some sum of
these angles and also verify the same inequality.  For any $n$ and
any $\varepsilon$, we've shown the existence of a configuration
$\CC_{\epsilon}$ of $n$ points of $\Hyp$ such that
$\angle_\mm(\CC_\varepsilon)\geq \frac{\pi}{n}-\varepsilon$, which
proves $\alpha_n(\Hyp)=\frac{\pi}{n}$.

\section{Concluding remarks}

One could ask the question of finding $\alpha_n(M)$ for other
types of metric spaces or manifolds. When $M=\R^d$, there are a
few immediate observations that one could make. Notice that for a
configuration $\CC$ of $n$ points in a Euclidean space, the
inequality $\alpha_\mm(\CC)\leq\frac{\pi}{3}$ always holds and the
equality $\alpha_\mm(\CC)=\frac{\pi}{3}$ is equivalent to having a
configuration of $n$ points all at an equal distance from one
another. For $\R^d$, this is possible if and only if $n\leq d+1$.
(To see this consider a generalized tetrahedron, or formally a
full-dimensional simplex.) It follows that
$\alpha_{n}(\R^d)=\frac{\pi}{3}$ if and only if $n\leq d+1$.
Another problem which seems similar in nature to finding
$\alpha_{n}(\R^3)$ is for $M=\Sp^2$. This problem is related to
the problem of finding, for $n$ points on the sphere, a sharp
upper bound on the minimum distance between them, which is a known
difficult problem.
\bibliographystyle{plain}
\def\cprime{$'$}

\end{document}